\documentclass{amsart}

\usepackage{amsmath,amssymb,amscd,amsthm}
\usepackage{multicol,longtable,tabularx}
\usepackage{enumerate}
\usepackage{latexsym}
\usepackage{float}
\usepackage{amstext,graphicx,color}
\usepackage[english]{babel}
\usepackage{ifthen}
\usepackage{verbatim}
\usepackage{calc}
\usepackage[vcentermath]{youngtab}

\newtheorem{theorem}{Theorem}[section]
\newtheorem{lemma}[theorem]{Lemma}
\newtheorem{proposition}[theorem]{Proposition}
\newtheorem{corollary}[theorem]{Corollary}

\newtheorem{main lemma}[theorem]{Main Lemma}
\newtheorem{conjecture}[theorem]{Conjecture}

\theoremstyle{definition}
\newtheorem{definition}[theorem]{Definition}
\newtheorem{example}[theorem]{Example}

\newtheorem{remark}[theorem]{Remark}

\newcommand{\zn}{\mathbb{Z}_m}
\newcommand{\fx}{f(x_1,\ldots,x_n)}
\newcommand{\N}{\mathbb{N}}
\newcommand{\Z}{\mathbb{Z}}
\newcommand{\F}{F\langle X|G\rangle}
\newcommand{\UT}{UT_m(F)}

\newcommand{\con}{f\equiv0}
\newcommand{\FF}{F\langle X\rangle}

\newcommand{\g}{\Gamma_{m_1,\ldots,m_r}^G}
\newcommand{\G}{\Gamma_n^G}

\author{Lucio Centrone}\thanks {Lucio Centrone is partially supported by FAPESP 2013/06752-4}
\address{IMECC, Universidade Estadual de Campinas, Rua S\'ergio Buarque de Holanda, 651, 13083-859 Campinas (SP) Brazil}\email{centrone@ime.unicamp.br}\author{Alessio Cirrito}\address{Dipartimento Di Matematica, Universit\`a Degli Studi Di Palermo, Palermo, via Archirafi 34, 90123}\email{alessio.cirrito@unipa.it}
\begin{document}
\title[$Y$-proper graded cocharacters of upper-triangular matrices]{$Y$-proper graded cocharacters of upper-triangular matrices of order $m$ graded by the $m$-tuple $\phi=(0,0,1,\ldots,m-2)$}
\begin{abstract} Let $F$ be a field of characteristic 0. We consider the algebra $UT_m(F)$ of upper-triangular matrices of order $m$ endowed with an elementary $\Z_m$-grading induced by the $m$-tuple $\phi=(0,0,1,\ldots,m-2)$, then we compute its $Y$-proper graded cocharacter sequence and we give the explicit formulas for the multiplicities in the case $m=2,3,4,5$.\end{abstract}\maketitle

\section{introduction}

Let $F$ be a field of characteristic 0. The algebra of upper-triangular matrices $UT_m(F)$ is a central object in the theory of PI-algebras satisfying a non-matrix polynomial identity. In fact, the polynomial identities of $UT_m(F)$ may serve as a measure of the complexity of the polynomial identities of finitely generated algebras with non-matrix identity in the same way as the polynomial identities of $M_m(F)$ measure the complexity of the identities of arbitrary PI-algebras.

In light of the work of Kemer (see \cite{kem1,kem2,kem3}) we have that the knowledge of graded identities and related graded structures such as graded cocharacters and graded codimensions of PI-algebras could be useful in order to understand better the polynomial identities of algebras. For example, the multiplicities of the graded cocharacter sequence of a PI-algebra provide an upper bound for the multiplicities of the ordinary cocharacter sequence.

Coming back to upper-triangular matrices, in \cite{val1} Valenti described all $\Z_2$-gradings over $UT_2(F)$ and many numerical characteristics related to these graded identities. Recall that a grading over the matrix algebra $M_m(F)$ is said to be elementary if the matrix units are homogeneous. In \cite{vaz1} Valenti and Zaicev described all the $G$-gradings over upper-triangular matrices, provided that $G$ is a finite abelian group and it turned out that such gradings are all isomorphic to elementary ones.

Among the various techniques used in order to manipulate polynomials, a special role is played by proper polynomials, expecially for graded identities. In the last year of the first half of the previous century Specht (see \cite{spe1}) started to use proper polynomials when dealing with polynomial identities of unitary algebras. In \cite{dre3} Drensky started
using them in a quantitative approach to the study of $T$-ideals of free algebras, developing an idea used
by Volichenko \cite{vol1}. When dealing with $G$-graded algebras one considers the free algebra $F\langle Y\cup Z\rangle$, where the $Y$'s are variables of degree $e=1_G$ and the $Z$'s are variables of homogeneous degree $g\neq e$. As long as proper polynomials are linear combinations of products of commutators on the generators, we say that a polynomial $f\in F\langle Y\cup Z\rangle$ is $Y$-proper if all the $y$'s occurring in $f$ appear in commutators on the generators only. It is well known (see, for instance, Lemma 1 Section 2 in \cite{dnd1}) that all graded polynomial identities of a superalgebra follow from the $Y$-proper ones. In \cite{dkv1} Di Vincenzo, Koshlukov and Valenti obtained a description of the $Y$-proper polynomials in the relatively-free $G$-graded algebra of $UT_m(F)$, where $G$ is a finite group. Notice that Koshlukov and Valenti in \cite{Kos}, found a multilinear basis of the relatively-free $\Z_m$-graded algebra of $UT_m(F)$.

Concerning the cocharacter sequence of the upper-triangular matrices, the explicit knowledge of multiplicities is well known only for $UT_2(F)$. Recently in \cite{bod1} Boumova and Drensky gave an easy algorithm which calculates the generating function of the cocharacter sequence of $UT_m(F)$. In \cite{cec1} Centrone and Cirrito gave a combinatorial method in order to compute the exact value of the $Y$-proper graded cocharacter sequence of $UT_m(F)$ endowed with the grading induced by that studied by Di Vincenzo in \cite{Divi05} and by Vasilovsky in \cite{Vas01} and they showed the exact value of the multiplicities for upper-triangular matrices of small size. We recall that the latter grading is an elementary grading over $UT_m(F)$ induced by the $m$-tuple $\psi=(0,1,\ldots,m-1)$ (see \cite{Divi05} and \cite{Vas01}).

Following this line of research, in this paper we deal again with elementary gradings on $UT_m(F)$. We are interested in studying the behaviour of the $Y$-proper cocharacter and codimension sequences of $UT_m(F)$ when the elementary $\Z_m$-grading is given by an $m$-tuple of the type $(\underbrace{0,\ldots,0}_{l},g_2,\ldots,g_{m-l+1})$. We observe that the grading of Vasilovsky is obtained for $l=1$ and all the other $g_i$'s are pairly different. Notice also that in the case $m=l$ the graded polynomial identies of $UT_m(F)$ are simply the ordinary polynomial identities of $UT_m(F)$. In particular, we consider the $\Z_m$-grading over $UT_m(F)$ induced by the $m$-tuple $\phi=(0,0,1,\ldots,m-2)$ and we ask what happens to the $Y$-proper graded cocharacter sequence. As well as in \cite{cec1} we give a combinatorial method in order to compute the $Y$-proper cocharacter sequence of $UT_m(F)$. The paper is organized as follows: the first section is introductory and furnishes all the tools for the study of graded polynomial identities and related structures. In the second section we focus on the grading induced by $\phi$ and we state our main theorem. The sections from 3 to 5 are dedicated to a computational inspection of the multiplicities of the $Y$-proper cocharacters for upper-triangular matrices of small size. We want to mention that there is a simple formula (in the spirit of Proposition 7.4 below) which allows to calculate the graded cocharacter sequence knowing the $Y$-proper graded cocharacters. We finish with a comparison between the experimental results obtained in \cite{cec1} and we observe that the multiplicities obtained when dealing with the grading induced by $\psi$ are asymptotically greater than in the case the grading is induced by $\phi$.

\section{Graded structures}

We introduce the key tools for the study of graded polynomial identities. We start off with the following definition:

\begin{definition} Let $G$ be a group and $A$ be an algebra over a field $F.$ We say that the algebra $A$ is \emph{$G$-graded} if $A$ can be written as the direct sum of subspaces $A=\bigoplus_{g\in G}A^{g}$ such that for all $g,\ h\in G$ one has $A^{g}A^{h}\subseteq A^{gh}.$\end{definition}

It is easy to note that, if $a$ is any element of $A$, it can be uniquely written as a finite sum $a=\sum_{g\in G}a_g,$ where $a_g\in A^{g}.$ We shall refer to the subspaces $A^{g}$ as \textit{$G$-homogeneous components} of $A$. An element $a\in A$ is said to be of \emph{homogeneous $G$-degree} $g$, written $\|a\|=g$, if $a\in A^{g}$. We always write $a^{g}$ if $a\in X^{g}$. Accordingly, an element $a\in A$ is called \emph{$G$-homogeneous} if exists $g\in G$ such that $a\in A^{g}.$ If $B\subseteq A$ is a subspace of $A,$ $B$ is $G$-graded if and only if $B=\bigoplus_{g\in G}(B\cap A^{g}).$ In the same way one defines \textit{$G$-graded} \emph{algebras, subalgebras, ideals, etc.}

Let $\{X^{g}\mid g \in G\}$ be a family of disjoint countable sets of indeterminates. We consider $X=\bigcup_{g\in G}X^{g}$ and we denote by $\F$ the free-associative algebra freely-generated by $X$. An indeterminate $x\in X$ is said to be of \emph{homogeneous $G$-degree} $g$, written $\|x\|=g$, if $x\in X^{g}$. We always write $x^{g}$ if $x\in X^{g}$. The homogeneous $G$-degree of a monomial $m=x_{i_1}x_{i_2}\cdots x_{i_k}$ is defined to be $\|m\|=\|x_{i_1}\|\cdot\| x_{i_2}\|\cdot\cdots\cdot\|x_{i_k}\|$. For every $g \in G$ we denote by $\F^g$ the subspace of $\F$ spanned by all the monomials having homogeneous $G$-degree $g$. Notice that $\F^g\F^{g'}\subseteq \F^{gg'}$ for all $g,g' \in G$. Thus \[\F=\bigoplus_{g\in G}\F^g\] is a $G$-graded algebra. We refer to the elements of the $G$-graded algebra $\F$ as \textit{$G$-graded polynomials} or, simply, \textit{graded polynomials}.

\begin{definition} If $A$ is a $G$-graded algebra, a $G$-graded polynomial $\fx$ is said to be a \emph{($G$-graded) polynomial identity} of $A$ if $$f(a_1,a_2,\cdots,a_n)=0$$ for all $a_1,a_2,\cdots,a_n\in\bigcup_{g\in G}A^g$ such that $a_k\in A^{\parallel x_k\parallel}$, $k=1,\cdots,n$. In this case we also say that $f\equiv0$ on $A$. Moreover $A$ is said to be a \textit{PI-graded algebra} if there exists a non-zero graded polynomial $f$ that is a graded polynomial identity for $A$.\end{definition}

For any $G$-graded algebra $A$ we consider \[T_G(A):=\{f\in\F|\text{\rm$\con$ on $A$}\},\]i.e., the set of $G$-graded polynomial identities of $A$.

\begin{definition} An ideal $I$ of $\F$ is said to be a \textit{$T_{G}$-ideal} if it is invariant under all $F$-endomorphisms $\varphi:\F\rightarrow\F$ such that $\varphi\left(\F^g\right)\subseteq\F^g$ for all $g\in G$.\end{definition}

Hence $T_G(A)$ is a $T_G$-ideal of $\F.$ On the other hand, it is easy to check that any $T_G$-ideal of $\F$ is $T_G(A)$ for some $G$-graded algebra $A$.

The theory of PI-graded algebras passes through the representation theory of the symmetric group. More precisely, it is useful to study the following $S_n$-modules.
\begin{definition} The elements of the set \[P_n^G=\textrm{span}\langle x_{\sigma(1)}^{g_1}x_{\sigma(2)}^{g_2}\cdots x_{\sigma(n)}^{g_n}|g_i\in G, \sigma\in S_n\rangle\]
are called \emph{$G$-graded multilinear polynomials} of degree $n$ of $\F$.\end{definition}

It turns out that $P_n^G$ is a left $S_n$-module under the natural left action of the symmetric group $S_n$; we denote the $S_n$-character of the factor module $P_n^G(A):=P_n^G/(P_n^G\cap T_G(A))$ by $\chi_n^G(A),$ and by $c_n^G(A)$ its dimension over $F$. We say: \[\left(\chi_n^G(A)\right)_{n\in\N}\ \textrm{is the $G$-\emph{graded cocharacter sequence} of $A$};\]
\[\left(c_n^G(A)\right)_{n\in\N}\ \textrm{is the $G$-\emph{graded codimension sequence} of $A$}.\]

If $G=\{g_1,\ldots,g_r\}$ is a finite group with unit $e$, for any $l_{g_1},\ldots,l_{g_r}\in \N$ let us consider the space of multilinear polynomials in the indeterminates labelled as follows: $x_1^{g_1},\ldots,x_{l_{g_1}}^{g_1},$ then $x_{l_{g_1}+1}^{g_2},\ldots,x_{l_{g_1}+l_{g_2}}^{g_2}$ and so on. We denote this linear space by $P_{l_{g_1},\ldots,l_{g_r}}^G$. The latter is a left $S_{l_{g_1}}\times\cdots\times S_{l_{g_r}}$-module. We shall denote by $\chi_{l_{g_1},\ldots,l_{g_r}}^G(A)$ the character of the left module $P_{l_{g_1},\ldots,l_{g_r}}^G(A)/(P_{l_{g_1},\ldots,l_{g_r}}^G(A)\cap T_G(A))$ and by $c_{l_{g_1},\ldots,l_{g_r}}^G(A)$ its dimension.

Since the ground field $F$ is infinite, a standard \emph{Vandermonde-argument} yields that a polynomial $f$ is a $G$-graded polynomial identity for $A$ if and only if its multihomogeneous components (with respect to the natural $\N$-grading), are identities as well. Moreover, since $\textrm{char}(F)=0,$ the well known multilinearization process shows that the $T_G$-ideal of a $G$-graded algebra $A$ is determined by its multilinear polynomials, i.e., $P_n^G(A)$ is determined by the $P_{l_{g_1},\ldots,l_{g_r}}^G(A)$. We remark that, given the character $\chi_{l_{g_1},\ldots,l_{g_r}}^G(A),$ the graded cocharacter $\chi_n^G(A)$ is known as well. More precisely, we have the following result by Di Vincenzo (see \cite{Divi07}, Theorem 2).

\begin{theorem}\label{gradedcocharacters}
Let $A$ be a PI $G$-graded algebra, then \[\chi_n^G(A)=\sum_{\begin{array}{c}
                      (l_{g_1},\ldots,l_{g_r}) \\
                      l_{g_1}+\ldots+l_{g_r}=n
                    \end{array}}\chi_{l_{g_1},\ldots,l_{g_r}}^G(A)^{\uparrow S_n}.\] Moreover \[c_n^G(A)=\sum_{\begin{array}{c}
                    (l_{g_1},\ldots,l_{g_r}) \\
                    l_{g_1}+\ldots+l_{g_r}=n
                    \end{array}}{n\choose l_{g_1},\ldots,l_{g_r}}c_{l_{g_1},\ldots,l_{g_r}}^G(A).\]\end{theorem}

Suppose we are dealing with a $G$-graded algebra and let us consider the free algebra $F\langle Y\cup Z\rangle$ (where $Y$ is the set of all indeterminates of $G$-degree $1_G$ and $Z$ is the set of all the remaining indeterminates). The $Y$-proper polynomials (see \cite{Dr06}, Section 2; \cite{Divi06}, Section 2) are the elements of the unitary $F$-subalgebra
$B$ of $\FF$ generated by the elements of $Z$ and by all non-trivial commutators. Roughly speaking, a
polynomial $f\in F\langle Y\cup Z\rangle$ is $Y$-proper if all the $y\in Y$ occurring in $f$ appear in commutators only. Notice
that if $f\in F\langle Z\rangle,$ then $f$ is $Y$-proper.
It is well known (see, for instance, Lemma 1 Section 2 in \cite{Divi06}) that all the graded polynomial identities
of a superalgebra $A$ follow from the $Y$-proper ones. This means that the set $T_{\Z_2}(A)\cap B$ generates the
whole $T_{\Z_2}(A)$ as a $T_{\Z_2}$-ideal. Similarly, for any group $G$, all the $G$-graded polynomial identities of a $G$-graded algebra $A$ follow from the $Y$-proper ones. This means that the set $T_G(A)\cap B$ generates the whole $T_G(A)$ as a $T_G$-ideal. Let us define $B(A):=B/(T_G(A)\cap B)$. We shall refer to $B(A)$ as \textit{$Y$-proper relatively-free algebra} of $A$.

We shall denote by $\Gamma_{n}^G$ the set of multilinear $Y$-proper polynomials of $P_{n}^G$. It is not difficult to see that $\G$ is a left $S_n$-submodule of $P_{n}^G$ and the same holds for $\G\cap T_G(A)$. Hence the factor module \[\G(A):={\G}/({\G\cap T_G(A)})\] is a $S_n$-submodule of $P_{n}^G(A)$; we denote the $S_n$-character of the factor module $\G/(\G\cap T_G(A))$ by $\xi_n^G(A)$, and by $\gamma_n^G(A)$ its dimension over $F$. We say: \[\left(\xi_n^G(A)\right)_{n\in\N}\ \textrm{is the $G$-\emph{graded proper cocharacter sequence} of $A$};\]
\[\left(\gamma_n^G(A)\right)_{n\in\N}\ \textrm{is the $G$-\emph{graded proper codimension sequence} of $A$}.\]
We shall denote by $\Gamma_{m_1,\ldots,m_r}^G$ the set of multilinear $Y$-proper polynomials of $P_{m_1,\ldots,m_r}^G$ such that $m=\sum_{i=1}^{r-1}m_i$. We observe that $\g$ is a left $S_{m_1}\times\cdots\times S_{m_r}$-submodule of $P_{m_1,\ldots,m_r}^G$ and the same holds for $\g\cap T_G(A)$. Hence the factor module \[\g(A):={\g}/({\g\cap T_G(A)})\] is a $S_{m_1}\times\cdots\times S_{m_r}$-submodule of $P_{m_1,\ldots,m_r}(A)^G$. We denote the $S_{m_1}\times\cdots\times S_{m_r}$-character of the factor module $\g/(\g\cap T_G(A))$ by $\xi_{m_1,\ldots,m_r}^G(A),$ and by $\gamma_{m_1,\ldots,m_r}(A)$ its dimension over $F$. When referring to $A$ without any ambiguity, we shall use $\gamma_{m_1,\ldots,m_r}$ instead of $\gamma_{m_1,\ldots,m_r}(A)$.

Following word by word the proof of Di Vincenzo in \cite{Divi07}, we have the analog of Theorem \ref{gradedcocharacters} for the graded-proper cocharacters and codimensions.

\begin{proposition}\label{prop2.1}
Let A be a PI G-graded algebra, then \[\xi_n^G(A)=\sum_{\begin{array}{c}
                      (l_{g_1},\ldots,l_{g_r}) \\
                      l_{g_1}+\ldots+l_{g_r}=n
                    \end{array}}\xi_{l_{g_1},\ldots,l_{g_r}}^G(A)^{\uparrow S_n}.\] Moreover \[\gamma_n^G(A)=\sum_{\begin{array}{c}
                    (l_{g_1},\ldots,l_{g_r}) \\
                    l_{g_1}+\ldots+l_{g_r}=n
                    \end{array}}{n\choose l_{g_1},\ldots,l_{g_r}}\gamma_{l_{g_1},\ldots,l_{g_r}}^G(A).\]\end{proposition}

We recall that a partition of the non-negative integer $n$ is a sequence of integers $\lambda=(\lambda_1,\ldots,\lambda_r)$ such that \[\text{\rm $\lambda_1\geq\cdots\geq\lambda_r>0$  and  $\lambda_1+\cdots+\lambda_r=n.$}\] In this case we shall write \[\lambda\vdash n.\] We assume two partitions $\lambda=(\lambda_1,\ldots,\lambda_r)$ and $\mu=(\mu_1,\ldots,\mu_s)$ to be equal if $r=s$ and \[\lambda_1=\mu_1,\ldots,\lambda_r=\mu_r.\]
When $\lambda=(\lambda_1,\ldots,\lambda_{k_1+\cdots+k_p})$ and \[\lambda_1=\cdots=\lambda_{k_1}=\mu_1,\ldots,\lambda_{k_1+\cdots+k_{p-1}+1}=\cdots=\lambda_{k_1+\cdots+k_p}=\mu_p,\] we accept the notation \[\lambda=(\mu_1^{k_1},\ldots,\mu_p^{k_p}).\]

\begin{definition}
Given a partition $\lambda=(\lambda_1,\ldots,\lambda_r)$, we associate to $\lambda$ the skew tableau $[\lambda]$ having $r$ rows such that its $i$-th row contains $\lambda_i$ squares. We say $[\lambda]$ is the \textit{Young diagram} of $\lambda$.
\end{definition}

\begin{example}If $\lambda=(5,3,1^2)$, then\[[\lambda]=\yng(5,3,1,1).\] In what follows we shall indicate by $\begin{tabular}{|l|l|l|l|}
\hline
\multicolumn{3}{|l|}{$k$}\\
\hline
\cline{1-2}
\end{tabular}$ the ``strip'' partition $(k)$.\end{example}

In what follows we use the following notation: if $G$ is a group and $F$ a field, we denote by $F[G]$, the $F$-group algebra of $G$.
By the Thorem of Maschke, if $G$ is a finite group, every finite dimensional representation is completely reducible, i.e., the group algebra $F[G]$ is semisimple and isomorphic to the direct sum of matrix algebras with entries from division algebras. Moreover, every finite dimensional left $G$-module is a direct sum of irreducible $G$-modules that are isomorphic to a minimal left ideal of $F[G]$. If $G=S_n$, the symmetric group of order $n$, the left irreducible $S_n$-modules (and their related characters) may be described in terms of partitions and Young diagrams. Indeed this applies for (ordinary, graded, Y-proper graded) cocharacters of PI-algebras. With abuse of notation we shall denote the irreducible character associated to the partition $\lambda$ by its Young diagram $[\lambda]$. When no ambiguity occurs we denote the irreducible character associated to $\lambda$ by $\lambda$ as well. We also recall that cocharacters of commutative algebras lie in ``strips''.

\section{Gradings on upper-triangular matrices: graded identities and cocharacters}

We give a general overlook at well known results about the graded identities of the algebra $UT_m(F)$ of upper-triangular matrices over a field $F$ graded by a finite abelian group. Then we focus on a particular $\Z_m$-grading obtained by the $m$-tuple $\phi=(0,0,1,2,\ldots,m-2)$.

Given an $m$-tuple $\overline{g}=(g_1,\ldots,g_m)$ of elements of a group $G$, the \textit{elementary $G$-grading} induced by $\overline{g}$ on the matrix algebra $M_m(F)$ is defined by $M_m(F)^{g_k}=span_F\langle e_{ij}|g_ig_j^{-1}=g_k\rangle$ for each $k\in\{1,\ldots,m\}$. If $A$ is a subalgebra of $M_m(F)$ which is $G$-graded with respect to this grading, we say that $\overline{g}$ is an elementary $G$-grading on $A$.

In what follows we shall recall some well known results about the graded polynomial identities of upper-triangular matrices endowed with any elementary $G$-grading, where $G$ is a finite abelian group. We start with the following definition.

\begin{definition}We shall call \emph{normal} any $Y$-commutator $c$ of $Z$-degree at most 1 such that $c=z$ or $c=[z,y_{i_1},\ldots,y_{i_t}]$ for some $t\geq1$.

Moreover we say:
\begin{itemize}
\item a normal commutator $[y_{i_1},\ldots,y_{i_t}]$ of $Z$-degree 0 is \emph{semistandard} if the indices $i_1,\ldots,i_p$ satisfy the inequalities $i_1>i_2\leq\cdots\leq i_p$.
\item a normal commutator $[z,y_{i_1},\ldots,y_{i_t}]$ of $Z$-degree 1 is \emph{semistandard} if the indices $i_1,\ldots,i_p$ satisfy the inequalities $i_1\leq i_2\leq\cdots\leq i_p$.
\end{itemize}\end{definition}

Let us consider $\UT$ endowed with an elementary $G$-grading. Consider the following definition (see \cite{Divi02}):

\begin{definition}
Let $\epsilon=(\epsilon_1,\ldots,\epsilon_m)$ be an elementary $G$-grading over $UT_m(F)$ and let $\overline{\eta}=(\eta_{1}, \ldots, \eta_{l})\in G^{l}$. We say $\overline{\eta}$ is a \textit{good sequence} (or \textit{$\epsilon$-good}) with respect to $\epsilon$ if there exists a sequence of $l$ matrix units $(r_{1}, \ldots, r_{l})$ in the Jacobson radical of $UT_{m}(F)$ such that the product $r_{1} \cdots r_{l}$ is not zero and $\eta_i$ is the homogeneous degree of $r_{i}$ for each $i=1, \ldots,l$. Otherwise $\overline{\eta}$ is called \textit{$\epsilon$-bad sequence}.
\end{definition}

For any sequence $\overline{\eta}\in G^{m}$ we consider the polynomial $f_{\overline{\eta}}=f_{\overline{\eta},1}f_{\overline{\eta},2}\cdots f_{\overline{\eta},m}$, where $f_{\overline{\eta},i}=[x_{e,2i-1},x_{e,2i}]$ if $\eta_{i}=e$ while $f_{\overline{\eta},i}=x_{\eta_{i}}$ if $\eta_{i}\neq e$. The following theorem was established in \cite{dkv1}. It holds over any infinite field.

\begin{theorem}\label{diViVa}
Let $G$ be a finite abelian group and $\epsilon=(\epsilon_{1}, \ldots, \epsilon_{m})$ be an elementary $G$-grading on $UT_{m}(F)$. Then the $T_G$-ideal $T_{G}(UT_{m}(F), \epsilon)$ of $G$-graded polynomial identities of $UT_{m}(F)$ is generated by all multilinear polynomials $f_{\overline{\eta}}$, where $\overline{\eta}=(\eta_{1}, \ldots, \eta_{l})$ lies on the set of all $\epsilon$-bad sequences and $l \leq m$.
\end{theorem}

In \cite{cec1} the authors studied the $Y$-proper graded cocharacters of $UT_m(F)$ equipped with the grading studied by Di Vincenzo and Vasilovsky. We recall that the latter grading over $UT_m(F)$ is a $\Z_m$-grading induced by the $m$-tuple $\psi=(0,1,\ldots,m-1)$. We want to study the case of $UT_m(F)$ graded by the $m$-tuple $\phi=(0,0,1,2,\ldots,m-2)$ of $\Z_m$. We start this investigation giving the complete description of $T_{\Z_d}(UT_d(F))$, $d=3,4,5$ as a Corollary of Theorem \ref{diViVa}.

\begin{corollary}
Let $G=\Z_m$ and let us denote by $y_{i}$ the variables of homogeneous degree $0$, by $z_{i}$ the variables of homogeneous degree $1$, by $t_{i}$ the variables of homogeneous degree $2$ and by $r_{i}$ the variables of homogeneous degree $3$. Then

\begin{enumerate}
  \item \[T_{G}(UT_{3}(F))=\langle [y_{1},y_{2}][y_{3},y_{4}],z_{1}z_{2},z[y_{1},y_{2}]\rangle^{T_G}.\]

  \item \[T_{G}(UT_{4}(F))=\langle [y_{1},y_{2}][y_{3},y_{4}],zt,[y_{1},y_{2}]z,t_{1}t_{2},t[y_{1},y_{2}],z_{1}z_{2}z_{3}\rangle^{T_G}.\]

  \item \[T_{G}(UT_{5}(F))=\langle [y_{1},y_{2}][y_{3},y_{4}],z[y_{1},y_{2}],zr,t[y_{1},y_{2}],t_{1}t_{2},tr,r[y_{1},y_{2}],\]\[
rz,rt,r_{1}r_{2},z_{1}z_{2}t,z_{1}z_{2}r,
      z_{1}tz_{2},tz_{1}z_{2},t_{1}zt_{2}\rangle^{T_G}.\]

\end{enumerate}
\end{corollary}
\proof
The proof is immediate by using Theorem \ref{diViVa}. More precisely, the possible bad sequences for $UT_{3}(F)$ are: $(0,0)$, $(0,1)$, $(1,1)$.
Analogously the bad sequences for $UT_{4}(F)$ are: $(0,0)$, $(0,1)$, $(1,2)$, $(2,2)$, $(0,2,0)$, $(1,1,1)$, $(2,0,2)$.
Finally the possible bad sequences for $UT_{5}(F)$ are: $(0,0)$, $(1,0)$, $(1,3)$, $(2,0)$, $(2,2)$, $(2,3)$, $(3,0)$, $(3,1)$, $(3,2)$, $(3,3)$, $(1,1,2)$, $(1,1,3)$, $(1,2,1)$, $(2,1,1)$, $(2,1,2)$.
\endproof

\begin{remark}\label{usethat}
The polynomial $z[y_1,y_2]$ is always a graded polynomial identity for $UT_m(F)$.
\end{remark}

The following result was proved for $UT_3(F)$ by Cirrito in \cite{cir1}. Following word by word his proof, we may generalize it for $UT_m(F)$, $m\geq3$.

\begin{lemma}\label{lem:14}
Let $UT_{m}(F)$ be equipped with the elementary $G$-grading induced by the $m$-tuple $\phi=(e,e,g_1,g_2,\ldots,g_{m-2})$ such that for each $i\neq j$ we have $g_i\neq g_j$ and $g_i\neq e$. If we denote by $ y_{i}$ the variables of homogeneous degree $e$ and by $z_{j}$ the variables of homogeneous degree $g\neq e$, then
\begin{enumerate}
  \item $[y_{i_{1}}, \ldots, y_{i_{h}},z,y_{j_{1}}, \ldots, x_{j_{k}}]\equiv \sum_{l \in I}\alpha_{l}g_{l} \quad \mbox{modulo}\quad T_{G}(UT_{m}(F))$
where $I$ is a finite set of indices, $h \geq 0$, $k \geq 2$, $\alpha_{l} \in F$ and $g_{l}$ is a product of two commutators for all $l \in I$. The first commutator is in the $y_{i}$'s only, and the second is a commutator in the $y_{i}$'s and the variable $z$ shifted in the first position (i.e.: $g_{l}=[y_{k_{1}}, \ldots, y_{k_{t}}][z,y_{r_{1}}, \ldots, y_{r_{m}}]$.

  \item $$[z,y_{i_{1}}, \ldots, y_{i_{k}},y_{2},y_{1},y_{j_{1}}, \ldots, y_{j_{h}}]\equiv$$ $$[z,y_{i_{1}}, \ldots, y_{i_{k}},y_{1},y_{2},y_{j_{1}}, \ldots, y_{j_{h}}]+\sum_{l \in I}\alpha_{l}g_{l} \quad \mbox{modulo}\quad T_{G}(UT_{m}(F)),$$
where $i_{1}, \ldots, i_{h}$ are not necessarily ordered indices, $h,k \geq 0$, $I$ is a finite set of indices, $\alpha_{l} \in F$ and $g_{l}$ is a product of two commutators. The first commutator is in the $y_{i}$'s only, and the second is a commutator in the $y_{i}$'s and $z$, for all $l \in I$.

  \item $$[z,y_{i_{1}}, \ldots, y_{i_{h}},\overline{y}_{1},y_{j_{1}}, \ldots, y_{j_{l}},\overline{y}_{2},y_{t_{1}}, \ldots, y_{t_{m}}]\equiv  \sum_{l \in I}\alpha_{l}g_{l}$$ $\mbox{modulo}\quad T_{G}(UT_{m}(F)),$
where $h,l,m \geq 0$, $I$ is a finite set of indices, $\alpha_{l} \in F$ for all $l \in I$ and $g_{l}$ is a product of two commutators with no alternating variables. Here $-$ means alternation on corresponding elements.
\end{enumerate}
\end{lemma}

In \cite{Divi02} Di Vincenzo, Koshlukov and Valenti obtained the following description of the $Y$-proper polynomials in the relatively-free graded algebra $\FF/(\FF\cap T_G(\UT))$.

\begin{theorem}\label{DiViKosVa}
A linear basis for the $Y$-proper polynomials in the relatively-free graded algebra $\FF/(\FF\cap T_G(\UT))$ consists of 1 and of the polynomials $c_1\cdots c_k$ where each polynomial $c_i$ is a semistandard commutator and the sequence $\tilde{\mu}_c=({\|c_1\|,\ldots,\|c_k\|})$ is good.
\end{theorem}

Notice that Koshlukov and Valenti in \cite{Kos}, found a multilinear basis of the relatively-free algebra $\FF/(\FF\cap T_{\zn}(\UT)).$ We have the following lemma that we shall use later on.

\begin{lemma}\label{properlemma} Let $m\geq2$ and consider $UT_m(F)$ with its $\phi$-grading. Then for any $n\in\N$ and for any $\sigma\in S_n,$ \[[z,y_{\sigma(1)},\ldots,y_{\sigma(n)}]\equiv_{T_{\zn}(UT_m(F))}[z,y_1,\ldots,y_n]+\sum_{i\in I}\alpha_ig_i,\]where $g_{i}$ is a product of two commutators for all $i \in I$, where the first commutator is semistandard of degree 0 and the second is a semistandard commutator of degree $\|z\|$.\end{lemma}

\proof
We proceed by induction on $n$. If $n=2$, then $[z,y_2,y_1]=[z,y_1,y_2]+z[y_2,y_1]-[y_2,y_1]z$ and we are done using Remark \ref{usethat}. Suppose the result true for $n-1>2$, then we prove it for $n$. By Lemma \ref{lem:14}, we have \[[z,y_{\sigma(1)},\ldots,y_{\sigma(n)}]=[z,y_1,\ldots,y_n]+\sum\alpha_lg_l,\] where the $g_l$'s are product of two commutators one of degree 0 and the other one with one $z$ as well as in Lemma \ref{lem:14} with length strictly less than that of $[z,y_1,\ldots,y_n]$. Now the lemma follows by induction and by the well known fact that the semistandard polynomials of degree 0 are a basis of the relatively-free proper algebra of $UT_2(F)$ (see for example \cite{dre4}).
\endproof

Given a sequence $\alpha=(\alpha_1,\ldots,\alpha_k)$ of elements of $\zn,$ we can associate to $\alpha$ the \textit{$m$-th uple of multiplicities} \[\mu(\alpha)=(\mu_1,\ldots,\mu_m)\] such that for any $i=1,\ldots,m$, \[\text{\rm $\mu_i=$ number of $\alpha_j$ such that $\alpha_j=[i-1]_m$}.\] We observe that distinct sequences are allowed to have the same $m$-th uple of multiplicities. Then \[\text{\rm $\alpha$ is a good sequence if and only if $\mu_1=0$ and $\sum_{j=2}^m\mu_j(j-1)\leq m-2$}.\] Fix now $\tilde{l}=(l_1,\ldots,l_m)$ such that $\sum_{j=2}^m l_j(j-1)\leq m-1$ and let \[S_{\tilde{l}}=\{\alpha=(\alpha_1,\ldots,\alpha_k)|\mu(\alpha)=(0,l_2,\ldots,l_m)\}.\]
 If $M$ is a $F[G]$-module and $m\in M$, we denote by $F[G]m$ the action of $F[G]$ on $m$. We also suggest that in what follows we implicitely use the fact that two representations afford the same $G$-character if they represent isomorphic $F[G]$-modules. With this in mind, we have the following.

\begin{theorem}\label{propcoch} Let $m\in\N$ and consider $UT_m(F)$ with its $\phi$-grading. Then for any $l_1,\ldots,l_m\in\N$ such that $\sum_{j=2}^ml_j(j-1)\leq m-2$,\[\xi_{l_1,\ldots,l_m}^{\zn}(UT_m(F))=\]\[\sum_{(\alpha_1,\ldots,\alpha_k)\in S_{\tilde{l}}}
\sum_{s_1+\ldots+s_{k+1}=l_1}
\left([\lambda(s_1)]\otimes\underbrace{\begin{tabular}{|c|cc|c|}
\hline
\empty & \multicolumn{2}{|c|}{$\cdots$} & \empty\\
\hline
\cline{1-2}
\end{tabular}}_{s_2}\otimes\cdots\otimes\underbrace{\begin{tabular}{|c|cc|c|}
\hline
\empty & \multicolumn{2}{|c|}{$\cdots$} & \empty\\
\hline
\cline{1-2}
\end{tabular}}_{s_{k+1}}\right)^{\uparrow S_{l_1}}\otimes$$$$\otimes(\yng(1)\otimes\cdots\otimes\yng(1))^{\uparrow S_{l_2}}\otimes\cdots\otimes(\yng(1)\otimes\cdots\otimes\yng(1))^{\uparrow S_{l_m}},\]where $\lambda(s_1)$ is the partition  $(s_1-1,1)$ with associated irreducible character \[\yng(5,1).\] \end{theorem}
\proof Let $(\alpha_1,\ldots,\alpha_k)\in S_{\tilde{l}}$ such that $\mu(\alpha)=(0,l_2,\ldots,l_m)$. Let $l_1\in\N$ such that $\sum_{i=1}^m l_i=n.$ Put $A=UT_m(F)$ and let $f\in\Gamma_{l_1,\ldots,l_m}\subseteq\Gamma_n^{\zn}(A)$. If $H=S_{l_1}\times\cdots\times S_{l_m},$ then $F[H]f^{\uparrow S_n}$ has dimension $|S_n:H|\dim_F F[H]f$. On the other hand, by Lemma \ref{properlemma}, $F[S_n]f$ is generated by \[[y_{l_1},\ldots,y_{{l_{s_1}}}][z^{\alpha_i},y_{i_1},\ldots,y_{{i_{s_2}}}]\cdots[z^{\alpha_j},y_{j_1},\ldots,y_{j_{{s_{k+1}}}}],\] where $l_1>l_2< l_3<\cdots<l_{s_1}$, the indices of the $y$'s of the other commutators are strictly increasing and $s_1+\cdots+s_{k+1}=l_1$. These commutators span the proper component for the ordinary polynomial identities of $UT_2(F)$ and the corresponding $S_{s_1}$-character is that for partition $(s_1-1,1)$. The latter polynomials are linearly independent, hence $\dim_FF[S_n]f=|S_n:H|\dim_F F[H]f$. Now the proof follows since the $S_n$-action in both of the situations is the same and the action on the $z$'s is the trivial one.\endproof

\section{Graded proper cocharacters of $UT_2(F)$ and $UT_3(F)$}
We use the results obtained in the previous section to compute the $Y$-proper graded cocharacter and codimension sequences for $UT_2(F)$ and $UT_3(F)$.

We start with the following combinatorial lemmas whose proofs repeat word by word those of Lemmas 4.1 and 4.2 of \cite{cec1}.

\begin{lemma}\label{calcolo1} Let $n\geq 1,$ then \[\sum_{l=1}^{n-1}([(l,1)]\otimes[(n-l-1)])^{\uparrow S_{n}}=\sum_{\begin{array}{c}
a+b+c=n\\
c\leq1
\end{array}}(a-b+1)[(a,b,c)].\]\end{lemma}

\begin{lemma}\label{calcolo2} Let $n\geq1,$ then \[\sum_{s=0}^{n-2}\sum_{t=1}^{n-s-1}([(s)]\otimes[(t,1)]\otimes[(n-s-t-1)])^{\uparrow S_{n}}\]\[=
\sum_{\begin{array}{c}
a+b+c+d=n\\
d\leq1
\end{array}}\sum_{i_1=b}^a\sum_{i_2=c}^b\sum_{i_3=d}^c(i_1-i_2+1)[(a,b,c,d)]=\sum_{\begin{array}{c}
a+b+c+d=n\\
d\leq1
\end{array}}C(a,b,c,d),\]where \[C=c(\frac{a-b+1}{2}((b+1)(a+2)+c(c-1))+\frac{c}{2}(b(b+1)-(a+1)(a+2)))\] if $d=1$ and \[C=(c+1)(\frac{a-b+1}{2}((b+1)(a+2)+c(c-1))+\frac{c}{2}(b(b+1)-(a+1)(a+2)))\] if $d=0$.
\end{lemma}

\begin{theorem}
Let $n\geq2$, then \[\xi_{n}^{\Z_2}(UT_2(F))=[(n-1,1)]\] and \[\gamma_n^{\Z_2}(UT_2(F))=n-1.\]
\end{theorem}
\proof Let $n\geq2$ and $G=\Z_2$, then we have no good sequences and the $Y$-proper $G$-graded cocharacters of $UT_2(F)$ are the $Y$-proper cocharacters of $UT_2(F)$ with the trivial grading and we are done. Analogously for the codimension sequence.\endproof

\begin{theorem}\label{upper3}
Let $n\geq2$, then \[\xi_{n}^{\Z_2}(UT_3(F))=2[(n-1,1)]+[(n-2,1^2)]+[(n-2,2)]\] and \[\gamma_n^{\Z_2}(UT_3(F))=2n-1.\]
\end{theorem}
\proof Let $n\geq2$ and $G=\Z_3$, then the only good sequence is $(1)$ that corresponds to the normal semistandard commutator $[z,y_1,\ldots,y_{n-1}]$, so $\gamma_{n-1,1,0}=1$ and by Proposition \ref{prop2.1} one has $\gamma_n^G(UT_3(F))=2n-1$.
Obviously, \[\xi_{n-1,1,0}^G(UT_3(F))=[(n-2,1)]\otimes[(1)]\] so by Theorem \ref{propcoch}, \[\xi_n^G(UT_3(F))=\left([(n-2,1)]\otimes[(1)]\right)^{\uparrow S_n}\otimes \emptyset =[(n-1,1)]+[(n-2,1^2)]+[(n-2,2)].\]Now the proof follows once we add $\xi_{n,0,0}^G(UT_3(F))$.\endproof

\section{Graded cocharacters of $UT_4(F)$}
We compute now the $Y$-proper graded cocharacters of $UT_4(F)$.

\begin{theorem}\label{mia1}
Let $n\geq2$, then $$\gamma_n^{\Z_4}(UT_4(F))=3n+2^{n-2}n(n-1)-1.$$ Moreover, $$\xi_{n}^{\Z_4}(UT_4(F))=\sum_{\lambda\vdash n}m_\lambda[\lambda],$$ where $$m_\lambda=\left\{\begin{array}{cc}
a+2 & if \lambda=(a,1), a \geq 4\\
6 & if \lambda=(3,2) \\
3(a-1) & if \lambda=(a,2), a \geq 4\\
5 & if \lambda=(b,b), b \geq 3 \\
8 & if \lambda=(b+1,b) , b \geq 3\\
4(a-b+1) & if \lambda = (a,b), b \geq 3, a \geq b+2\\
3a-2 & if \lambda=(a,1,1), a \geq 2 \\
5 & if \lambda=(2,2,1) \\
12 & if \lambda=(3,2,1) \\
7a-9 & if \lambda=(a,2,1), a > 3 \\
9 & if \lambda=(b,b,1), b >3 \\
15 & if \lambda=(b+1,b,1), b >3 \\
8(a-b+1) & if \lambda=(a,b,1), b \geq 3, a \geq b+2 \\
5 & if \lambda=(b,b,2), b \geq 2 \\
5(a-b+1) & if \lambda=(a,b,2), a \geq b+2\\
1 & if \lambda=(1^4)\\
3a-2 & if \lambda=(a,1^{3}), a \geq 2 \\
5 & if \lambda=(b,b, 1^{2}) \\
5(a-b+1) & if \lambda=(a,b,1^2), a \geq b+1, b \geq 2\\
2(a-b+1) & if \lambda=(a,b,2,1), a \geq b+1, b \geq 2\\
a-b+1 & if \lambda=(a,b,3), a \geq b+1, b \geq 3\\
\end{array}\right.$$
\end{theorem}
\proof
Let $n\geq2$ and $G=\Z_4$, then the good sequences are $(1),(2)$ and $(1,1)$ that correspond respectively to the normal semistandard commutators \[[z^{1},y_1,\ldots,y_{n-1}],\] \[[z^{2},y_1,\ldots,y_{n-1}]\] and \[[z_{\sigma(1)}^{1},y_{i_1},\ldots,y_{i_k}][z_{\sigma(2)}^{1},y_{j_1},\ldots,y_{j_{l}}],\ \ \ \sigma\in S_2,\] where $\{i_1,\ldots,i_k\}\cup\{j_1,\ldots,j_l\}=\{1,2,\ldots,n-2\}$ and ${k+l=n-2}$. It is easy to see that $\gamma_{n-1,1,0,0}^G=\gamma_{n-1,0,1,0}^G=1.$ Standard combinatorial facts allow us to say that $\gamma_{n-2,2,0,0}^G=2\sum_{i=0}^{n-2}{n-2\choose i}=2\cdot2^{n-2}=2^{n-1}$. By Proposition \ref{prop2.1}, we have $$\gamma_n^G(UT_4(F))=$$$$={n\choose n,0,0,0}(n-1)+{n\choose n-1,1,0,0}+{n\choose n-1,0,1,0}+2^{n-1}{n\choose n-2,2,0,0}=$$$$=3n+2^{n-2}n(n-1)-1.$$

Obviously, $$\xi_{n-1,1,0,0}^G(UT_4(F))=(n-1,1)\otimes(1)\otimes\emptyset \otimes \emptyset,$$ $$\xi_{n-1,0,1,0}^G(UT_4(F))=(n-1,1)\otimes\emptyset\otimes(1)\otimes \emptyset$$ and  $$\xi_{n-2,2,0,0}^G(UT_4(F))=\sum_{l=1}^{n-3}((l,1)\otimes(n-3-l))^{\uparrow S_{n-2}}\otimes((2)+(1^2))\otimes\emptyset \otimes \emptyset.$$ By Lemma \ref{calcolo1} we have:
\begin{equation}\label{coccolo}
W=\sum_{l=0}^{n-3}((l,1)\otimes(n-3-l))^{\uparrow S_{n-2}}=\sum_{l=1}^{n-3}((l,1)\otimes(n-3-l))^{\uparrow S_{n}}=\sum_{\begin{array}{c}
a+b+c=n-2\\
c\leq1
\end{array}}(a-b+1)(a,b,c).
\end{equation}

Let $$W_1=\left(W\otimes(2)\right)^{\uparrow S_n}$$ and $$W_2=\left(W\otimes(1^2)\right)^{\uparrow S_n},$$ then $$\xi_{n-2,2,0,0}^G(UT_4(F))=W_1+W_2.$$ Let us consider $W_1$. Firstly, we note that $W_1=\sum_{\lambda\vdash n}m_\lambda\lambda,$ where $\lambda$ is a partition of height at most 4. We observe that $\lambda=(n)$ for any $n \geq 1$ is not allowed as a partition of $W_1$. If $\lambda=(a,1)$, $a \geq4$ it comes from $(a-2,1)\otimes(2)$ with multiplicity 1, then due to Equation (\ref{coccolo}), $\lambda$ has total multiplicity $a-2$. We observe that $\lambda=(1^2), (2,1)$ are not allowed as a partitions of $W_1$. If $\lambda=(3,2),$ it comes from $(2,1)\otimes(2)$, with multiplicity 1, then due to Equation (\ref{coccolo}), $\lambda$ has total multiplicity $2$, whereas if $\lambda=(a,2),$ it comes from $(a-1,1)\otimes(2)$ or $(a-2,2)\otimes(2)$, with multiplicity 1, then due to Equation (\ref{coccolo}), it has total multiplicity $a-1+a-3=2a-4$. If $\lambda=(b,b)$, $b \geq 3$, it comes from $(b,b-2)\otimes(2)$ or $(b-1,b-1)\otimes(2)$, with multiplicity 1, then due to Equation (\ref{coccolo}), $\lambda$ has total multiplicity $3+1=4$ and if $\lambda=(b+1,b)$, $b\geq 3$, it comes from $(b+1,b-2)\otimes(2)$ or $(b,b-1)\otimes(2)$, with multiplicity 1, then due to Equation (\ref{coccolo}), it has total multiplicity $4+2=6$.
If $\lambda=(a,b)$, $b \geq 3$, $a \geq b+2$, it comes from $(a,b-2)\otimes(2)$, $(a-1,b-1)\otimes(2),$ or $(a-2,b)\otimes(2)$, with multiplicity 1, then due to Equation (\ref{coccolo}), $\lambda$ has total multiplicity $a-b+3+a-b+1+a-b-1=3(a-b+1)$. We observe that $\lambda=(1^3)$ is not allowed as a partition of $W_1$. If $\lambda=(a,1^2)$, $a\geq 3$ it comes from $(a-1,1)\otimes(2)$ or $(a-2,1,1) \otimes (2)$ with multiplicity 1, then due to Equation (\ref{coccolo}), $\lambda$ has total multiplicity $a-1+a-2=2a-3$.
We use the same techniques to obtain multiplicities for partitions $\lambda= (2,2,1), (3,2,1)$.
If $\lambda=(a,2,1)$, $a >3$, it comes from $(a,1) \otimes (2)$ or $(a-1,2)\otimes(2)$ or $(a-2,2,1)\otimes(2)$ or $(a-1,1,1) \otimes (2)$, with multiplicity 1, then due to Equation (\ref{coccolo}), $\lambda$ has total multiplicity $a+a-2+a-3+a-1=4a-6$. If we consider partitions of the type $\lambda=(b,b,1)$, $b>3$, they come from $(b,b-1)\otimes(2)$ or $(b-1,b-1,1)\otimes(2)$ or $(b,b-2,1) \otimes (2)$, with multiplicity 1, then due to Equation (\ref{coccolo}), they have total multiplicity $2+1+3=6$; whereas the partitions $\lambda=(b+1,b,1)$, $b>3$, come from $(b+1,b-1,1)\otimes(2)$ or $(b,b)\otimes(2)$ or $(b,b-1,1) \otimes (2)$ or $(b,b-2,1) \otimes (2)$, with multiplicity 1, then due to Equation (\ref{coccolo}), they have total multiplicity $3+1+2+3=9$.
If $\lambda=(a,b,1)$, $b> 3$, $a \geq b+2$, then it comes from $(a,b-1)\otimes(2)$ or $(a-1,b)\otimes(2)$ or $(a-1,b-1,1) \otimes (2)$ or $(a-2,b,1) \otimes (2)$ or $(a,b-2,1) \otimes (2)$, each with multiplicity 1, then due to Equation (\ref{coccolo}), $\lambda$ has total multiplicity $a-b+2+a-b+a-b+1+a-b-1+a-b+3=5(a-b+1)$.
If $\lambda=(b,b,2)$, $b \geq 2$, it comes from $(b,b)\otimes(2)$ or $(b,b-1,1)\otimes(2)$, with multiplicity 1, then due to Equation (\ref{coccolo}), $\lambda$ has total multiplicity $1+2=3$.
If $\lambda=(a,b,2)$, $b \geq 2$ $a \geq b+1$, it comes from $(a,b)\otimes(2)$ or $(a-1,b,1)\otimes(2)$ or $(a,b-1,1) \otimes (2)$, with multiplicity 1, then due to Equation (\ref{coccolo}), $\lambda$ has total multiplicity $a-b+1+a-b+a-b+2=3(a-b+1)$.
If $\lambda=(b,b,1,1)$, $b \geq 2$, it comes from $(b,b-1,1)\otimes(2)$, with multiplicity 1, then due to Equation (\ref{coccolo}), $\lambda$ has total multiplicity $2$.
If $\lambda=(a,b,1,1)$, $b \geq 2$, $a \geq b+1$ it comes from $(a-1,b,1)\otimes(2)$ or $(a,b-1,1)\otimes(2)$, with multiplicity 1, then due to Equation (\ref{coccolo}), $\lambda$ has total multiplicity $a-b+a-b+2=2(a-b+1)$.
Finally we can computes multiplicities for partitions $\lambda= (a,b,2,1)$ and $\lambda = (a,b,3)$ observing that they can be obtained from partitions $(a,b,1) \otimes (2)$ and $(a,b,1) \otimes (2)$ respectively.

Now we compute $W_2$. We note that $W_2=\sum_{\lambda\vdash n}m_\lambda\lambda,$ where $\lambda$ is a partition of height at most 5. We observe that $\lambda=(a), (a,1)$ for any $a \geq 1$, are not allowed as a partitions of $W_2$. If $\lambda=(a,b)$, $b \geq 2$, it comes from $(a-1,b-1)\otimes(1^2)$ with multiplicity 1, then due to Equation (\ref{coccolo}), $\lambda$ has total multiplicity $a-b+1$. If $\lambda=(a,b,1)$, it comes from $(a,b-1)\otimes(1^2)$ or $(a-1,b)\otimes(1^2)$ or $(a-1,b-1,1)\otimes(1^2)$ , with multiplicity 1, then due to Equation (\ref{coccolo}), $\lambda$ has total multiplicity $a-b+2+a-b+a-b+1=3(a-b+1)$.
If $\lambda=(a,1,1)$, it comes from $(a-1,1)\otimes(1^2)$, with multiplicity 1, then due to Equation (\ref{coccolo}), $\lambda$ has total multiplicity $a-1$.
If $\lambda=(b,b,1)$, $b \geq 2$, it comes from $(b,b-1)\otimes(1^2)$ or $(b-1,b-1,1)\otimes(1^2)$ , with multiplicity 1, then due to Equation (\ref{coccolo}), $\lambda$ has total multiplicity $2+1=3$.
If $\lambda=(b,b,2)$, $b \geq 2$, it comes from $(b,b-1,1)\otimes(1^2)$, with multiplicity 1, then due to Equation (\ref{coccolo}), $\lambda$ has total multiplicity $2$.
If $\lambda=(a,b,2)$, $b \geq 2$, $a \geq b+1$, it comes from $(a-1,b,1)\otimes(1^2)$ or $(a,b-1,1)\otimes(1^2)$, with multiplicity 1, then due to Equation (\ref{coccolo}), $\lambda$ has total multiplicity $a-b+a-b+2=2(a-b+1)$.
Note that partitions $\lambda=( a,b, c )$, $c \geq 3$, are not allowed as partitions of $W_2$.
If $\lambda=(1^4)$, it comes from $(1,1)\otimes(1^2)$, with multiplicity 1, then due to Equation (\ref{coccolo}), $\lambda$ has total multiplicity $1$.
If $\lambda=(b,b,1,1)$, $b \geq 2$, it comes from $(b,b)\otimes(1^2)$ or $(b,b-1,1)\otimes(1^2)$, with multiplicity 1, then due to Equation (\ref{coccolo}), $\lambda$ has total multiplicity $1+2=3$.
If $\lambda=(a,1^3)$, it comes from $(a,1)\otimes(1^2)$ or $(a-1,1,1)\otimes(1^2)$, with multiplicity 1, then due to Equation (\ref{coccolo}), $\lambda$ has total multiplicity $a+a-1=2a-1$.
If $\lambda=(a,b,1^2)$, $a \geq b+1$, $b \geq 2$ it comes from $(a,b)\otimes(1^2)$ or $(a,b-1,1)\otimes(1^2)$ or $(a-1,b,1)\otimes(1^2)$ with multiplicity 1, then due to Equation (\ref{coccolo}), $\lambda$ has total multiplicity $a-b+1+a-b+2+a-b=3(a-b+1)$.
Finally if $\lambda=(a,b,2,1)$, it comes from $(a,b,1)\otimes(1^2)$, with multiplicity 1, then due to Equation (\ref{coccolo}), $\lambda$ has total multiplicity $a-b+1$.

Now we summarize the previous partial results. If $W=\sum_{\lambda\vdash n}m_\lambda[\lambda],$ then $$m_\lambda=\left\{\begin{array}{cc}
a-2 & if \lambda=(a,1), a \geq 4\\
4 & if \lambda=(3,2) \\
3a-5 & if \lambda=(a,2), a \geq 4\\
5 & if \lambda=(b,b), b \geq 3 \\
8 & if \lambda=(b+1,b) , b \geq 3\\
4(a-b+1) & if \lambda = (a,b), b \geq 3, a \geq b+2\\
3a-4 & if \lambda=(a,1,1), a \geq 2 \\
5 & if \lambda=(2,2,1) \\
12 & if \lambda=(3,2,1) \\
7a-9 & if \lambda=(a,2,1), a > 3 \\
9 & if \lambda=(b,b,1), b >3 \\
15 & if \lambda=(b+1,b,1), b >3 \\
8(a-b+1) & if \lambda=(a,b,1), b \geq 3, a \geq b+2 \\
5 & if \lambda=(b,b,2), b \geq 2 \\
5(a-b+1) & if \lambda=(a,b,2), a \geq b+2\\
1 & if \lambda=(1^4)\\
3a-2 & if \lambda=(a,1^{3}), a \geq 2 \\
5 & if \lambda=(b,b, 1^{2}) \\
5(a-b+1) & if \lambda=(a,b,1^2), a \geq b+1\\
2(a-b+1) & if \lambda=(a,b,2,1), a \geq b+1\\
a-b+1 & if \lambda=(a,b,3), a \geq b+1\\
\end{array}\right.$$ Now the proof follows once we add $\xi_{n-1,1,0,0}^G(UT_4(F))$, $\xi_{n-1,0,1,0}^G(UT_4(F))$ and $\xi_{n,0,0,0}^G(UT_4(F))$.
\endproof
\section{graded cocharacters of $UT_5(F)$}
We consider now $UT_5(F)$. In this case we calculate the exact value of the multiplicities for partitions $(a,b)$ and $(a,1^{b})$, $a \geq 1$, $b \geq 0$.
Instead, we will suggest the algorithm for partitions $(a,b,c)$, $a\geq b \geq c \geq 0$, $(a,b,c,1)$, $a\geq b \geq c \geq 1$, $(a,b,c,2)$, $a\geq b \geq c \geq 2$, $(a,b,c,3)$, $a\geq b \geq c \geq 3$,
$(a,b,c,2,1)$, $a\geq b \geq c \geq 2$.

\begin{theorem}\label{mia2}
Let $n\geq3$, then \[\gamma_n^{\Z_5}(UT_5(F))=4n+n(n-1)+32^{n-3}n(n-1)+3^{n-3}n(n-1)(n-2)-1.\]
Moreover, if we consider partitions $(a,1^{b})$, $a\geq 1$ and $b \geq 0$, then $$\xi_{n}^{\Z_5}(UT_5(F))=\sum_{\lambda\vdash n}m_\lambda[\lambda],$$
where $$m_\lambda=\left\{\begin{array}{cc}
a^{2}-2a+11 & \text{\rm if $\lambda=(a,1)$ and $a\geq3$}\\
23          &  if \lambda=(3, 1^{2})\\
3a^{2}-3a+5 & \text{\rm if $\lambda=(a,1^{2})$ and $a\geq4$}\\
17          &  if \lambda=(2, 1^{3})\\
34          &  if \lambda=(3, 1^{3})\\
4a^{2}-a & \text{\rm if $\lambda=(a,1^{3})$ and $a\geq4$}\\
3 & if \lambda=(1^5)\\
12          &  if \lambda=(2, 1^{4})\\
3a^{2}-2a+2 & \text{\rm if $\lambda=(a,1^{4})$ and $a\geq3$}\\
2          &  if \lambda=( 1^{6})\\
2a^{2} & \text{\rm if $\lambda=(a,1^{5})$ and $a\geq2$}\\
\frac{a(a+1)}{2} & \text{\rm if $\lambda=(a,1^{6})$ and $a\geq1$}.\\
\end{array}\right.$$
\end{theorem}

\proof
Let $n\geq3$ and $G=\Z_5$, then the good sequences are $(1)$, $(2)$, $(3)$, $(1,1)$, $(1,2)$ and $(1,1,1)$ that correspond respectively to the normal semistandard commutators
\[[z^{1},y_1,\ldots,y_{n-1}],\] \[[z^{2},y_1,\ldots,y_{n-1}],\] \[[z^{3},y_1,\ldots,y_{n-1}]\] \[[z_{\sigma(1)}^{1},y_{i_1},\ldots,y_{i_k}][z_{\sigma(2)}^{1},y_{j_1},\ldots,y_{j_{l}}],\ \sigma\in S_2,\]
\[[z^{1},y_{i_1},\ldots,y_{i_k}][z^{2},y_{j_1},\ldots,y_{j_{l}}],\ \ [z^{2},y_{i_1},\ldots,y_{i_k}][z^{1},y_{j_1},\ldots,y_{j_{l}}]\]
\[[z_{\sigma(1)}^{1},y_{i_1},\ldots,y_{i_k}][z_{\sigma(2)}^{1},y_{j_1},\ldots,y_{j_{l}}][z_{\sigma(3)}^{1},y_{v_1},\ldots,y_{v_{m}}],\ \sigma\in S_3,\]
where \[\{i_1,\ldots,i_k\}\cup\{j_1,\ldots,j_l\}=\{1,2,\ldots,n-2\},\ {k+l=n-2}\] and \[\{i_1,\ldots,i_k\}\cup\{j_1,\ldots,j_l\}\cup\{v_1,\ldots,v_m\}=\{1,2,\ldots,n-3\},\ {k+l+m=n-2}.\]
We have \[\gamma_{n-1,1,0,0,0}^G=\gamma_{n-1,0,1,0,0}^G=\gamma_{n-1,0,0,1,0}^G=1\] and \[\gamma_{n-2,2,0,0,0}^G=\gamma_{n-2,1,1,0,0}^G=2^{n-1}\]
Standard combinatorial facts allow us to say \[\gamma_{n-3,3,0,0,0}^G=6\sum_{i=0}^{n-3}\sum_{j=0}^{n-3-i}{n-3-i \choose j}{n-3\choose i}=2\cdot3^{n-2}.\]
By Proposition \ref{prop2.1} we have \[\gamma_n^G(UT_5(F))=4n+2^{n-1}n(n-1)+3^{n-3}n(n-1)(n-2)-1.\]
Now we focus on cocharacters. Obviously, $$\xi_{n-1,1,0,0,0}^G(UT_5(F))=(n-1)\otimes(1)\otimes\emptyset\otimes\emptyset \otimes\emptyset,$$ $$\xi_{n-1,0,1,0,0}^G(UT_5(F))=(n-1)\otimes\emptyset\otimes(1)\otimes\emptyset \otimes\emptyset,$$
$$\xi_{n-1,0,0,1,0}^G(UT_5(F))=(n-1)\otimes\emptyset\otimes\emptyset\otimes(1)\otimes\emptyset.$$

$\xi_{n-2,2,0,0,0}^G(UT_5(F))$ and $\xi_{n-2,1,1,0,0}^G(UT_5(F))$ have the same decomposition; by Theorem \ref{mia1} we obtain the following decomposition.

$$\xi_{n-3,3,0,0,0}^G(UT_5(F))=\sum_{s=0}^{n-2}\sum_{t=1}^{n-s-1}((s)\otimes(t,1)\otimes(n-s-t-1))^{\uparrow S_{n}}\otimes((2)+(1^3)+(2,1))\otimes\emptyset\otimes\emptyset.$$

By Lemma \ref{calcolo2} we have that \begin{equation}\label{coccolo2}W=\sum_{s=0}^{n-2}\sum_{t=1}^{n-s-1}((s)\otimes(t,1)\otimes(n-s-t-1))^{\uparrow S_{n}}=
\sum_{\begin{array}{c}
a+b+c+d=n\\
d\leq1
\end{array}}\sum_{i_1=b}^a\sum_{i_2=c}^b\sum_{i_3=d}^c(i_1-i_2+1)(a,b,c,d).\end{equation}

Let \[W_1=\left(W\otimes(3)\right)^{\uparrow S_n},\] \[W_2=\left(W\otimes(1^3)\right)^{\uparrow S_n},\] and \[W_3=\left(W\otimes(2,1)\right)^{\uparrow S_n},\] then \[\xi_{n-3,3,0,0,0}^G(UT_5(F))=W_1+W_2+W_3.\]

Let us consider $W_1$. Firstly, we note that $W_1=\sum_{\lambda\vdash n}m_\lambda\lambda,$ where $\lambda$ is a partition of height at most 5. We observe that $\lambda=(a)$ for any $a \geq 1$, $\lambda= (1, 1), (2,1), (3,1)$ are not allowed as partitions of $W_1$. If $\lambda=(a,1),$ $a \geq 4$, it comes from $(a-3,1)\otimes(3)$ with multiplicity 1, then due to Equation (\ref{coccolo2}), $\lambda$ has total multiplicity $(a-3)(a-1)$.
We observe that $\lambda=(1^{3}), (2,1,1)$ are not allowed as partitions of $W_1$. If $\lambda=(3,1,1)$, it comes from $(1,1)\otimes(3)$ with multiplicity 1, then due to Equation (\ref{coccolo2}), $\lambda$ has total multiplicity $3$. If $\lambda=(a,1,1)$, $a \geq 4$, it comes from $(a-3,1,1)\otimes(3)$ or $(a-2,1)\otimes (3)$, with multiplicity 1, then due to Equation (\ref{coccolo2}), $\lambda$ has total multiplicity $2(a-3)(a-1)+2-(a-2)(a-1)+a(a-2)$.
We observe that $\lambda=(1^{4}), (2,1^{3})$ are not allowed as partitions of $W_1$. If $\lambda=(3,1^{3})$, it comes from $(1,1,1)\otimes(3)$ with multiplicity $1$, then due to Equation (\ref{coccolo2}), $\lambda$ has total multiplicity $2$. If $\lambda=(a,1^{3}),$ $a \geq 4$, it comes from $(a-3,1^{3})\otimes(3)$ or $(a-2,1,1)\otimes (3)$, with multiplicity 1, then due to Equation (\ref{coccolo2}), $\lambda$ has total multiplicity $(a-3)(a-1)+1-(a-2)(a-1)+2a(a-2)+2-a(a-1)$.
In the same way $\lambda=(1^{5}), (2,1^{4})$ are not allowed as partitions of $W_{1}$.
If $\lambda=(a,1^{4}),$ $a \geq 3$, it comes from $(a-2,1^{3})\otimes(3)$, with multiplicity 1, then it has total multiplicity $a(a-2)+1-\frac{a(a-1)}{2}$.

Now we compute $W_{2}$. We note that $W_2=\sum_{\lambda\vdash n}m_\lambda\lambda,$ where $\lambda$ is a partition of height at most $6$.
We observe that $\lambda=(a,1)$, $a \geq 1$, $\lambda=(a,1,1)$, $a \geq 1$, $\lambda=(1^{4})$ are not allowed as a partitions of $W_{2}$.
If $\lambda=(a,1^{3}),$ $a \geq 2$, it comes from $(a-1,1)\otimes(1^{3})$, with multiplicity 1, then due to Equation (\ref{coccolo2}),
$\lambda$ has total multiplicity $(a+1)(a-1)$.
In the same way $\lambda=(1^{5})$ has total multiplicity $3$.
If $\lambda=(a,1^{4}),$ $a \geq 2$, it comes from $(a-1,1,1)\otimes(1^{3})$ or $(a,1) \otimes (1^{3})$, with multiplicity 1, then due to Equation (\ref{coccolo2}),
$\lambda$ has total multiplicity $a(a+2)+2(a-1)(a+1)+2-a(a+1)$.
If $\lambda=(1^{6})$, it comes from $(1^{3})\otimes(1^{3})$, with multiplicity 1 and then the total multiplicity is $2$.
If $\lambda=(a,1^{5}),$ $a \geq 2$, it comes from $(a-1,1^{3})\otimes(1^{3})$ or $(a,1,1) \otimes (1^{3})$, with multiplicity 1, then $\lambda$ has total multiplicity $2a(a+2)+2-(a+1)(a+2)+(a-1)(a+1)+1-\frac{a(a+1)}{2}$.
If $\lambda=(a,1^{6}),$ $a \geq 1$, it comes from $(a,1^{3})\otimes(1^{3})$, with multiplicity 1, then due to Equation (\ref{coccolo2}),
$\lambda$ has total multiplicity $\frac{a(a+1)}{2}$.

Now compute $W_{3}$. We note that $W_3=\sum_{\lambda\vdash n}m_\lambda\lambda,$ where $\lambda$ is a partition of height at most $5$.
We observe that $\lambda=(a)$, $a \geq 1$ is not allowed as a partition of $W_{3}$; also $\lambda=(a,1),$ $a \geq 1$ and $\lambda =(1^{3}), (2,1^{2})$ are not allowed as partitions of $W_{2}$.
If $\lambda=(3,1^{2})$, it comes from $(1,1)\otimes(2,1)$, with total multiplicity $3$.
If $\lambda=(a,1^{2}),$ $a \geq 3$, it comes from $(a-2,1)\otimes(2,1)$ then due to Equation (\ref{coccolo2}),
$\lambda$ has total multiplicity $a(a-2)$.
$\lambda =(1^{4})$ is not allowed as a partition of $W_{3}$.
If $\lambda=(2,1^{3})$, it comes from $(1^{2})\otimes(2,1)$, with total multiplicity $3$.
If $\lambda=(a,1^{3}),$ $a \geq 3$, it comes from $(a-2,1^{2})\otimes(2,1)$ or $(a-1,1) \otimes (2,1)$, with multiplicity 1, then due to Equation (\ref{coccolo2}),
$\lambda$ has total multiplicity $(a+1)(a-1)+2(a-1)a+2-a(a-1)$.
$\lambda=(1^{5})$ is not allowed as a partition of $W_{3}$.
If $\lambda=(2,1^{4})$, it comes from $(1^{3})\otimes(2,1)$, with total multiplicity $2$.
If $\lambda=(a,1^{4}),$ $a \geq 3$, it comes from $(a-1,1^{2})\otimes(2,1)$ and $(a-2, 1^{3}) \otimes (2,1)$, with multiplicity 1, then due to Equation (\ref{coccolo2}),
$\lambda$ has total multiplicity $a(a-2)+1-a(a-1)+2(a-1)(a+1)+2-a(a+1)$.
If $\lambda=(a,1^{5}),$ $a \geq 2$, it comes from $(a-1,1^{3})\otimes(2,1)$  with multiplicity 1, then due to Equation (\ref{coccolo2}),
$\lambda$ has total multiplicity $(a-1)(a+1)+1-\frac{a(a+1)}{2}$.

Now we summarize the previous partial results:
$$m_\lambda=\left\{\begin{array}{cc}
(a-3)(a-1) & \text{\rm if $\lambda=(a,1)$ and $a\geq3$}\\
6          &  if \lambda=(3, 1^{2})\\
3a^{2}-9a+6 & \text{\rm if $\lambda=(a,1^{2})$ and $a\geq4$}\\
6          &  if \lambda=(2, 1^{3})\\
20          &  if \lambda=(3, 1^{3})\\
4a^{2}-7a+4 & \text{\rm if $\lambda=(a,1^{3})$ and $a\geq4$}\\
3 & if \lambda =(1^5)\\
12          &  if \lambda=(2, 1^{4})\\
3a^{2}-2a+2 & \text{\rm if $\lambda=(a,1^{4})$ and $a\geq3$}\\
2          &  if \lambda=( 1^{6})\\
2a^{2} & \text{\rm if $\lambda=(a,1^{5})$ and $a\geq2$}\\
\frac{a(a+1)}{2} & \text{\rm if $\lambda=(a,1^{6})$ and $a\geq1$}.\\
\end{array}\right.$$

Now the proof follows once we add
\begin{gather*}
\xi_{n,0,0,0,0}^G(UT_5(F)),\\
\xi_{n-1,1,0,0,0}^G(UT_5(F)), \\
\xi_{n-1,0,1,0,0}^G(UT_5(F)), \\
\xi_{n-1,0,0,1,0}^G(UT_5(F)), \\
\xi_{n-2,2,0,0,0}^G(UT_5(F)), \\
\xi_{n-2,1,1,0,0}^G(UT_5(F)).
\end{gather*}
\endproof

\begin{theorem}
Let $n\geq3$, if we consider partitions $(a,b)$, $a,b\geq 1$, then \[\xi_{n}^{\Z_5}(UT_5(F))=\sum_{\lambda\vdash n}m_\lambda[\lambda],\]
where $$m_\lambda=\left\{\begin{array}{cc}
10 & if \lambda=(3,2)\\
37 & if \lambda=(4,2)\\
\frac{7a^{2}-11a+6}{2} & \text{\rm if $\lambda=(a,2)$ and $a\geq5$}\\
18 & if \lambda=(3,3)\\
52 & if \lambda=(4,3)\\
102 & if \lambda=(5,3)\\
7a^{2}-13a-6 & \text{\rm if $\lambda=(a,3)$ and $a\geq6$}\\
3a^{2}+1 & \text{\rm if $\lambda=(a,a)$ and $a\geq4$}\\
6a^{2}-6a-5 & \text{\rm if $\lambda=(a,a-1)$ and $a\geq5$}\\
3(3a^{2}-6a-4) & \text{\rm if $\lambda=(a,a-2)$ and $a\geq6$}\\
\frac{5a^{2}b-5ab^{2}+12ab-3a^{2}-3b^{2}+2b-2a-2}{2} & \text{\rm if $\lambda=(a,b)$ and $a,b\geq3$, $a-b \geq 3$}.\\
\end{array}\right.$$
\end{theorem}

\proof
As well as in the proof of Theorem \ref{mia2} we have to compute \[\xi_{n-3,3,0,0,0}^G(UT_5(F))=W_{1}+W_{2}+W_{3}\]
The multiplicities for partitions $(a)$, $a \geq 2$ and $(a,1)$, $a \geq 2$ are the same of those computed in the previous theorem.

We consider firstly $W_{1}$, then $W_1=\sum_{\lambda\vdash n}m_\lambda\lambda,$ where $\lambda$ is a partition of height at most $4$.
We observe that $\lambda=(2,2), (3,2)$ are not allowed as partitions of $W_{1}$.
If $\lambda=(4,2)$ it comes from $(2,1) \otimes (3)$ with multiplicity $1$, then due to Equation (\ref{coccolo2}), it has total multiplicity $8$.
If $\lambda=(a,2)$, $a \geq 5$, it comes from $(a-3,2) \otimes (3)$ or $(a-2,1) \otimes (3)$, with multiplicity $1$, then due to Equation (\ref{coccolo2}),
 it has total multiplicity $\frac{3}{2}(a-1)(a-4)+a(a-2)$.
If $\lambda=(4,3)$ it comes from $(3,1) \otimes (3)$ with multiplicity $1$, then due to Equation (\ref{coccolo2}), it has total multiplicity $15$.
If $\lambda=(5,3)$ it comes from $(3,2) \otimes (3)$ or $(4,1) \otimes (3)$ with multiplicity $1$, then due to Equation (\ref{coccolo2}), it has total multiplicity $39$.
If $\lambda=(a,3)$, $a \geq 6$, it comes from $(a-3,3) \otimes (3)$ or $(a-2,2) \otimes (3)$ or $(a-1,1) \otimes (3)$, with multiplicity $1$, then due to Equation (\ref{coccolo2}),
 it has total multiplicity $2(a-1)(a-5)+\frac{3}{2}a(a-3)+(a-1)(a+1)$.
If $\lambda=(a,a)$, $a \geq 4$, it comes from $(a,a-3) \otimes (3)$, with multiplicity $1$, then due to Equation (\ref{coccolo2}),
 it has total multiplicity $2(a-2)(a+2)$.
If $\lambda=(a,a-1)$, $a \geq 5$, it comes from $(a,a-4) \otimes (3)$ or $(a-1,a-3) \otimes (3)$, with multiplicity $1$, then due to Equation (\ref{coccolo2}),
 it has total multiplicity $\frac{5}{2}(a-3)(a+2)+\frac{3}{2}(a-2)(a+1)$.
If $\lambda=(a,a-2)$, $a \geq 6$, it comes from $(a,a-5) \otimes (3)$ or $(a-1,a-4) \otimes (3)$ or $(a-2,a-3) \otimes (3)$, with multiplicity $1$, then due to Equation (\ref{coccolo2}),
 it has total multiplicity $3(a-4)(a+2)+2(a-3)(a+1)+a(a-2)$.
If $\lambda=(a,b)$, $a,b \geq 3$, $a-b \geq 3$, it comes from $(a,b-3) \otimes (3)$ or $(a-1,b-2) \otimes (3)$ or $(a-2,b-1) \otimes (3)$ or $(a-3,b)\otimes (3)$, with multiplicity $1$, then due to Equation (\ref{coccolo2}),
 it has total multiplicity $\frac{(a+2)(b-2)}{2}(a-b+4)+\frac{(a+1)(b-1)}{2}(a-b+2)+\frac{ab(a-b)}{2}+\frac{(a-1)(b+1)}{2}(a-b-2)$.

Now we consider $W_{2}=(W \otimes (1^{3}))^{\uparrow S_{n}}$. We note that $\lambda=(a,b)$, for every $a \geq 1$, $b \geq 0$ is not allowed as a partition of $W_{2}$.
So we can consider $W_{3}=(W \otimes (2,1))^{\uparrow S_{n}}$.
Note that $\lambda=(2,2)$ it is not allowed as partition of $W_{3}$.
If $\lambda=(3,2)$ it comes from $(1,1) \otimes (2,1)$, with multiplicity $1$, then due to Equation (\ref{coccolo2}),
 it has total multiplicity $3$.
If $\lambda=(4,2)$ it comes from $(2,1) \otimes (2,1)$, with multiplicity $1$, then due to Equation (\ref{coccolo2}),
 it has total multiplicity $8$.
If $\lambda=(a,2)$, $a \geq 5$, it comes from $(a-2,1) \otimes (2,1)$, with multiplicity $1$, then due to Equation (\ref{coccolo2}),
 it has total multiplicity $a(a-2)$.
If $\lambda=(3,3)$ it comes from $(2,1) \otimes (2,1)$, with multiplicity $1$, then due to Equation (\ref{coccolo2}),
 it has total multiplicity $8$.
If $\lambda=(4,3)$ it comes from $(2,2) \otimes (2,1)$ or $(3,1) \otimes (2,1)$, with multiplicity $1$, then due to Equation (\ref{coccolo2}),
 it has total multiplicity $21$.
If $\lambda=(5,3)$ it comes from $(3,2) \otimes (2,1)$ or $(4,1) \otimes (2,1)$, with multiplicity $1$, then due to Equation (\ref{coccolo2}),
 it has total multiplicity $39$.
If $\lambda=(a,3)$, $a \geq 6$, it comes from $(a-2,2) \otimes (2,1)$ or $(a-1,1) \otimes (2,1)$, with multiplicity $1$, then due to Equation (\ref{coccolo2}),
 it has total multiplicity $\frac{3}{2}a(a-3)+(a+1)(a-1)$.
If $\lambda=(a,a)$, $a \geq 4$, it comes from $(a-1,a-2) \otimes (2,1)$, with multiplicity $1$, then due to Equation (\ref{coccolo2}),
 it has total multiplicity $(a+1)(a-1)$.
If $\lambda=(a,a-1)$, $a \geq 4$, it comes from $(a-1,a-3) \otimes (2,1)$ or $(a-2,a-2) \otimes (2,1)$, with multiplicity $1$, then due to Equation (\ref{coccolo2}),
 it has total multiplicity $\frac{3}{2}(a+1)(a-2)+\frac{a(a-1)}{2}$.
If $\lambda=(a,a-2)$, $a \geq 5$, it comes from $(a-1,a-4) \otimes (2,1)$ or $(a-2,a-3) \otimes (2,1)$, with multiplicity $1$, then due to Equation (\ref{coccolo2}),
 it has total multiplicity $2(a+1)(a-3)+a(a-2)$.
If $\lambda=(a,b)$, $a,b \geq 3$, $a-b \geq 3$, it comes from $(a-1,b-2) \otimes (2,1)$ or $(a-2,b-1) \otimes (2,1)$, with multiplicity $1$, then due to Equation (\ref{coccolo2}),
 it has total multiplicity $\frac{(a+1)(b-1)}{2}(a-b+2)+\frac{ab(a-b)}{2}$.

Now the proof follows if we sum these partial results and once we add
\begin{gather*}
\xi_{n,0,0,0,0}^G(UT_5(F)),  \\
\xi_{n-1,1,0,0,0}^G(UT_5(F)), \\
\xi_{n-1,0,1,0,0}^G(UT_5(F)), \\
\xi_{n-1,0,0,1,0}^G(UT_5(F)), \\
\xi_{n-2,2,0,0,0}^G(UT_5(F)), \\
\xi_{n-2,1,1,0,0}^G(UT_5(F)).
\end{gather*}
\endproof

\section{Conclusions}
We will draw out a deep difference between two gradings over the same algebra of the same group. We focus on the following definition.

\begin{definition} Let $\varphi$ and $\chi$ be characters of the symmetric group $S_n$. We consider the decompositions of $\varphi$ and $\chi$ into irreducible characters, say \[\varphi=\sum_{\lambda\vdash n}m_\lambda[\lambda],\ \ \ \chi=\sum_{\lambda\vdash n}m_\lambda'[\lambda]\] where $m_\lambda$ and $m_\lambda'$ are non-negative integers. We shall write \[\text{\rm $\varphi\leq\chi$ if and only if $m_\lambda\leq m_\lambda'$ for all $\lambda\vdash n$}.\]\end{definition}

It is well known (see \cite{Divi07} Proposition 1) that:

\begin{theorem}\label{minchar} Let $A$ be a $G$-graded algebra with ordinary cocharacter sequence $\chi_n(A)$. Then \[\chi_n(A)\leq\chi_n^G(A)\]\end{theorem}

Using similar arguments we can state the following:

\begin{theorem}\label{minchar2} Let $A$ be a $G$-graded algebra with proper cocharacter sequence $\zeta_n(A)$ and with Y-proper graded cocharacter sequence $\xi_n^G(A)$. Then \[\zeta_n(A)\leq\xi_n^G(A).\]\end{theorem}

We consider now the following result of Drensky (see for example \cite{dre4} Theorem 12.5.4).

\begin{proposition}\label{drensky1}
Let $A$ be a PI-algebra and $\chi_n(A)=\sum_{\lambda\vdash n}m_{\lambda}(A)[{\lambda}]$ its $n$-th cocharacter. Let $\xi_p(A)=\sum_{\nu\vdash p}k_{\nu}(A)[{\nu}]$ its p-th proper cocharacter, then $$m_{\lambda}(A)=\sum_{\nu\in S}k_{\nu}(A),$$ where $S=\{\nu=(\nu_1,\ldots,\nu_n)\mid\lambda_1\geq\nu_1\geq\lambda_2\geq\nu_2\geq\cdots\geq\lambda_n\geq\nu_n\}.$
\end{proposition}

As well as in \cite{cec1} notice that combining the computational results of the previous sections and Theorem \ref{minchar2} we are able, in principle, to provide an upper bound for the ordinary cocharacter sequence of $UT_m(F)$, at least for special partitions.

\begin{example}\label{ale}
Let $A:=UT_3(F)$ with its ordinary cocharacter sequence $\chi_n(A)=\sum_{\lambda\vdash n}m_\lambda[\lambda]$. We consider the partition $\lambda=(n-2,1^2)$. We use Proposition \ref{drensky1}, then $\lambda$ comes from $(n-2,1^2)\otimes1$, $(n-3,1^2)\otimes(1)$,$\ldots$, $(1^3)\otimes(n-3)$ or from $(n-2,1)\otimes(1)$, $(n-3,1)\otimes(2)$,$\ldots$, $(1^2)\otimes(n-2)$ where $(n-2,1^2)$, $(n-3,1^2)$,$\ldots$, $(1^3)$ and $(n-2,1)$, $(n-3,1)$,$\ldots$, $(1^2)$ represent irreducible modules of $\Gamma_n(A)$. By Theorem \ref{upper3} the total multiplicity is: \[3\sum_{k=1}^{n-2}1=3n-6.\]Then \[m_\lambda\leq3n-6.\]
\end{example}

\begin{example}
Let $A:=UT_4(F)$ with its ordinary cocharacter sequence $\chi_n(A)=\sum_{\lambda\vdash n}m_\lambda[\lambda]$. We consider the partition $\lambda=(n-2,1^2)$. We use Proposition \ref{drensky1}, then $\lambda$ comes from $(n-2,1^2)\otimes1$, $(n-3,1^2)\otimes(1)$,$\ldots$, $(2,1^2)\otimes(n-4)$ or from $(n-1,1)\otimes1$, $(n-2,1)\otimes(1)$,$\ldots$, $(2,1)\otimes(n-3)$ where $(n-2,1^2)$, $(n-3,1^2)$,$\ldots$, $(2,1^2)$ and $(n-1,1)$, $(n-2,1)$,$\ldots$, $(2,1)$ represent irreducible modules of $\Gamma_n(A)$. By Theorem \ref{mia1} the total multiplicity is: \[\sum_{k=1}^{n-2}3k-2+\sum_{k=4}^{n-1}k+2=2n^2-5n-7.\]Then \[m_\lambda\leq2n^2-5n-7.\]
\end{example}

We may compare these results with Examples 6.5 and 6.6 of \cite{cec1}. In both of the two cases, the multiplicities of the two partitions considered when dealing with the $\phi$-grading are asymptotically less than the multiplicities of the same partitions when dealing with the $\psi$-grading. This behavior may also be observed when dealing with growth functions of polynomial identities. In a forthcoming paper Aljadeff and Centrone proved that the Gelfand-Kirillov dimension of the $\Z_n$-graded relatively-free algebra of $M_n(F)$ is the greatest possible among all the $\Z_n$-gradings if equipped with the $\psi$-grading. For more details about the graded Gelfand-Kirillov dimension of graded algebras we refer to \cite{cen1}, \cite{cen2} and \cite{cen3}. In light of the previous facts we state the following conjecture.

\begin{conjecture}
Let $\chi^{\psi}_n=\sum m_\lambda^{\psi}[\lambda]$ be the $\Z_m$-graded cocharacter sequence of $UT_m(F)$ with the $\Z_m$-grading induced by the $\psi$-grading and $\chi_n=\sum m_\lambda[\lambda]$ be the $\Z_m$-graded cocharacter sequence of $UT_m(F)$ with any $\Z_m$-grading. Then for each $\lambda\vdash n$ we have $m_\lambda\leq m^{\psi}_\lambda$.
\end{conjecture}

\end{document}